\documentclass[preprint,10pt,leqno]{elsarticle}

\usepackage{fullpage}

\usepackage{amsrefs}

\usepackage{graphicx}
\usepackage{caption}
\usepackage{subcaption}
\usepackage{amssymb}
\usepackage{amsmath}
\usepackage{epstopdf}

\usepackage{amsthm}

\newcommand{\eps}{\epsilon}
\theoremstyle{plain}
\newtheorem{theorem}{Theorem}[section]
\newtheorem{theorem*}{Theorem}

\newproof{pf}{Proof}

\journal{}

\begin{document}

\begin{frontmatter}

\title{The Effects of Latent Infection on the Dynamics of HIV}


\author{Stephen Pankavich\fnref{support}}
\fntext[support]{The author was supported in part by NSF grants DMS-0908413 and DMS-1211667}

\address{Department of Applied Mathematics and Statistics, Colorado School of Mines, Golden, CO 80401}

\ead{pankavic@mines.edu}
           
\date{Received: date / Accepted: date}

\begin{abstract}
One way in which the human immunodeficiency virus (HIV-1) replicates within a host is by infecting activated CD$4+$ T-cells, which then produce additional
copies of the virus.  
Even with the introduction of antiretroviral drug therapy, which has been very successful over the past decade,  a large obstacle to the complete eradication of the virus is the presence of viral reservoirs in the form of latently infected CD$4+$ T-cells. 
We consider a model of HIV infection that describes T-cell and viral interactions, as well as, the production and activation of latently infected T-cells.  
Upon determining equilibrium states of the latent cell model, the local and global asymptotic behavior of solutions is examined, and the basic reproduction number of the system is computed to be strictly less than that of the corresponding three-component model, which omits the effects of latent infection.
In particular, this implies that a wider variety of parameter values will lead to viral eradication as $t \to \infty$ due to the appearance of latent CD$4+$ T-cells.  
With this realization we discuss possible alternative notions for eradication and persistence of infection other than traditional tools.
These results are further illustrated by a number of numerical simulations.

\end{abstract}

\begin{keyword}
HIV-1 \sep Mathematical model \sep Latently infected T-cells \sep Antiretroviral therapy \sep Global asymptotic stability
\end{keyword}

\end{frontmatter}

\section{Introduction}

The majority of cells infected with the human immunodeficiency virus (HIV-1) are activated CD$4+$ T-cells.
Once infected, these cells produce additional copies of virus, thereby prolonging the infection.
Upon detecting such an infection, the immune system mounts a complex adaptive response, controlling the virus population to a limited extent. 
Further control is available in the form of antiretroviral drugs, such as
Reverse Transcriptase Inhibitors (RTIs) and Protease Inhibitors (PIs) \cite{Janeway}. 
If such drugs are taken with sufficient frequency, the virus population is largely limited and
remains below the level of detection \cite{CF}. However, antiretroviral
therapy (ART) cannot fully eradicate the virus, as viral rebound occurs once therapy is interrupted
\cite{Arlen, Doyle}
A number of factors have been proposed to explain this viral rebound.
Most notably, it has been suggested that HIV lay dormant within a number of reservoirs.
Primary among these reservoirs are latently infected CD4+ T-cells \cite{Blankson}.
Though latently infected T-cells exist in the body with a much lower frequency than susceptible CD4+ T-cells,
the reservoir appears to decay very slowly, with a half-life measured between $6$ and $48$ months \cite{Ramratnam}.
Although infected, these cells do not produce new virions until activated, thus potentially providing a
longer-lived hiding place where the virus may evade control by either the immune system or antiretroviral
treatment \cite{Blankson}.
Consequently, in this paper, we analyze a mathematical model that includes latent infection and examine
the control of infection by ART.
We also assume that such latent T-cells exist significantly longer than productively infected CD4+ T-cells. 
Ultimately, we will show that a mathematical analysis of the most basic latent model demonstrates
that the inclusion of such cells increases the likelihood for viral clearance under the traditional approach of analyzing the basic reproduction number and the
associated stability of equilibria.  
While this will seem intuitive from a modeling perspective, as described later, it also appears contradictory to the widely-held notion that latently infected T-cells
are an important mechanism for the inability of ART to eradicate an established infection.
What this really implies is that standard mathematical tools are insufficient to realistically describe the dynamics of HIV when latent cells
are considered.
Instead, one must focus on the rate of decay of the infection, which is significantly slowed by the latent T-cell population.

A number of authors have studied the biological aspects of mathematical models concerning HIV dynamics that consider latently infected cells. 
%
Notably, Callaway and Perelson \cite{CallPerelson} studied low-level viremia,
Rong and Perelson \cite{RongPerelson} modeled viral blips and showed that
a latent reservoir could produce viral transients when activated by infection,
while Sedaghat et al. \cite{Sedaghat} employed a simple model for the dynamics of the latent
reservoir to show that its stability was unlikely to stem from ongoing replication during ART. 
In each of these studies, a reduced or linearized mathematical analysis was performed, but the nonlinear behavior of the associated model was not fully elucidated.
In the current study, we describe latently infected cells using a separate compartment, as did these authors, by assuming that a proportion of newly-infected cells become latently infected upon contact with the virus, but that they are not productively infected until they leave the latent state, which occurs at a rate $\alpha$ proportionate
to the strength of the latent cell population.
We note that the effects of viral mutation, which may continuously change model and parameter values, and the possible spatial dependence of parameters are ignored.
Using this model, we study the influence of the latent reservoir on the persistence of HIV infection and viral rebound.
Our results provide a new perspective on the methods of mathematical and stability analysis for viral and latent reservoir persistence.

The paper proceeds as follows.  In the next section, we will review some known results concerning the standard
three-component model of HIV dynamics.  In Section $3$, we introduce an additional population representing latently infected 
CD$4+$ T-cells, and study the effects that these cells have on the structure and behavior of the long-time dynamics of the model.  
In Section $4$, we discuss the ramifications of our results and, in particular, the need to construct more precise notions of viral eradication
and persistence.
The fifth section contains proofs of the theorems contained within previous sections.
In the final section, we conclude with a discussion of our results.

\section{The Three-component model}

In general, the modeling of HIV dynamics  \emph{in vivo} is
complicated by the appearance of spatial inhomogeneities, which can arise
from various reservoirs, such as those occurring within lymphatic
tissues \cite{KPAIDS, PH}.
Even when such inhomogeneities are ignored, however, 
these systems are often described to a sufficient degree by systems of ordinary differential equations
that include no spatial dependence. 
We begin by considering a three-component model for the evolution
of within-host HIV, that does not include spatial fluctuations or effects due to long-lived infected or latently infected cells.
This model has been widely-accepted as a descriptive representation
for the basic dynamics of HIV \cite{BCN, PKdB, TuckShip}.
It represents the populations of three components in a fixed volume at a given time $t$:
$T(t)$, the number of CD$4+$ T-cells that are susceptible to HIV-1
infection, $I(t)$ the number of infected T-cells that are actively producing virus particles,  and $V(t)$ the number of free virions. 
These quantities approximately satisfy the system of ordinary differential equations
\begin{equation}
\label{3CM}
\left \{ \begin{aligned}
\frac{dT} {dt} &= \lambda - d_T T - k T V \\
\frac{dI} {dt} &=  k T V - d_I I \\
\frac {dV} {dt} &= Nd_I I - d_V V.
\end{aligned} \right.
\end{equation}
Here $\lambda$ is the recruitment rate of susceptible T-cells and $d_T$ is their
mortality rate. The constant $k$ represents the rate of infection, which is included
within a bilinear mass action term, while $d_I$ is 
the death rate of productively infected cells and $d_V$ is the clearance rate of free virus. The parameter
$N$ is the burst size, i.e. the total number of virions produced by an infected cell during its life span.

\begin{figure*}[t]
\begin{center}
\includegraphics[scale = 0.58]{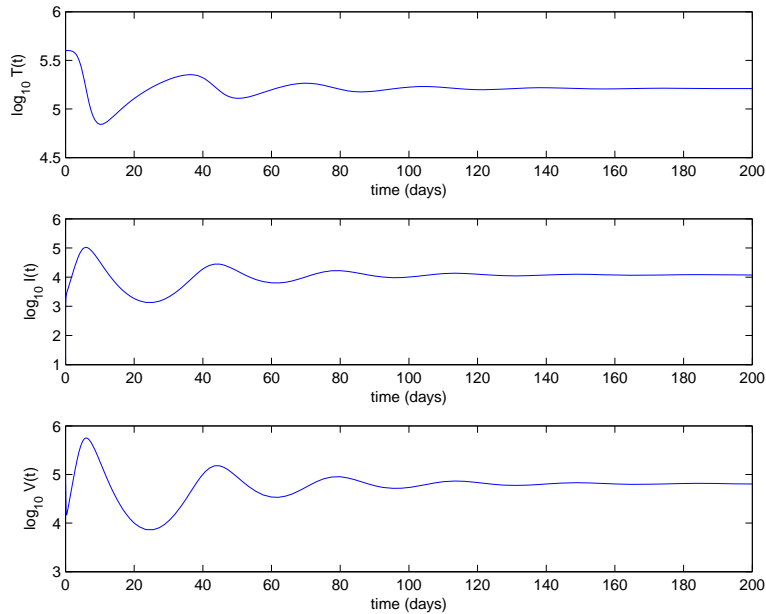}
\vspace{-0.2in}
\caption{A representative solution of (\ref{3CM}) with parameter values stated in Table \ref{tab:1}. 
The initial T-cell population is $T(0) = 4 \times 10^5$, while the initial viral load is $V(0) = 10^5$,
and $I(0) = 0$.
In this example, the system tends to the endemic equilibrium as $t \to \infty$ because $R_0 = 2.087$.
 }
\label{fig:1}  
\end{center}
\end{figure*}

The behavior of solutions to these equations has previously been analyzed in great detail \cite{SCCDHP, BCN, TuckWan, TuckShip}.
In particular, it is known that exactly two steady states exist, which we will write in the form $(T,I,V)$, namely a non-infective equilibrium 
$$E_{NI}: \qquad \left ( \frac{\lambda}{d_T}, 0, 0 \right )$$
and an infective or endemic equilibrium
$$E_I: \qquad \left ( \frac{\lambda}{d_T R_0}, \frac{d_T d_V}{kNd_I} (R_0 - 1) , \frac{d_T}{k} (R_0 - 1) \right )$$
where
$$R_ 0 = \frac{\lambda k N}{d_V d_T}.$$
The stability properties of these steady states are also well-known and depend only upon the single parameter $R_0$, called the basic reproduction number.
In particular, one can study the linearized analogue of (\ref{3CM}) and prove the local asymptotic stability of $E_{NI}$ if $R_0 \leq 1$ and
the local asymptotic stability of $E_I$ if $R_0 > 1$ \cite{TuckWan}.  This result effectively states that for starting values of the
populations which are close enough to the given equilibria, the solutions will tend to the respective equilibrium point
as $t \to \infty$.  Additionally, the global asymptotic stability of these equilibria is known.
In \cite{Korob} it was shown that initial populations are irrelevant in determining the long term dynamics of the solution.
More specifically, if $R_0 \leq 1$, then for any initial population of uninfected cells,
infected cells, and virions the solutions of (\ref{3CM}) tend to $E_{NI}$ as $t \to \infty$.  Contrastingly, if $R_0 > 1$ the same result holds for the
endemic equilibrium $E_I$. 
Figure \ref{fig:1} displays a representative graph of solutions for which $R_0 > 1$ and hence viral infection persists. 

\section{Inclusion of Latently-infected cells}

Though (\ref{3CM}) describes the basic mechanisms which account for the spread of HIV, it lacks the ability to describe the latent stage of a specific subpopulation of infected T-cells. 
Many studies \cite{Chun1995, Chun1997a, Chun1997b} have determined that upon infection and transcription of viral RNA into cell DNA, a fraction of CD$4+$ T-cells fail to actively produce virus until they are activated, possibly years after their initial infection.  
Such cells may possess a much longer lifespan than their counterparts, and are termed latently infected. 
Upon activation, latently infected cells do become actively productive, and hence begin to increase the viral load through
viral replication. 
A basic model of latent cell activation was initially developed to examine cell populations that contribute to the viral decline 
that occurs after administration of antiretroviral therapy \cite{PECVH}.
However, within \cite{PECVH} and other articles by related authors \cite{RongPerelson,RongPerelsonRev1, PKdB}, the mathematical analysis of the model is performed under a number of limiting assumptions, including a constant background population of susceptible T-cells and perfect efficacy of anti-retroviral drugs.
Thus, we focus on rigorously proving the resulting nonlinear dynamics without these assumptions.

As for (\ref{3CM}) we consider a model describing T-cells that may be susceptible or infected.  In addition, we let $L(t)$ represent the new population of latently infected T-cells that cannot produce virions at time $t$ but begin to do so once they are activated by recall antigens.  With this addition, the previously described three-component model now contains four components and is given by
\begin{equation}
\label{4CM}
\left \{ \begin{aligned}
\frac{dT} {dt} &= \lambda - d_T T - k T V \\
\frac{dI} {dt} &=  (1-p)k T V + \alpha L - d_I I \\
\frac{dL} {dt} &=  pk T V  - \alpha L - d_L L \\
\frac {dV} {dt} &= Nd_I I - d_V V.
\end{aligned} \right.
\end{equation}
Here, $p \in (0,1)$ is the proportion of infections that lead to the production of a latently infected T-cell, rather than a productively infected T-cell, and  $\alpha$ is the rate at which latently infected cells transition to become actively productive.
Additionally, $d_L$ is the rate at which latent cells are cleared from the system.
\begin{figure*}[t]
\begin{center}
\includegraphics[width = 5.8in]{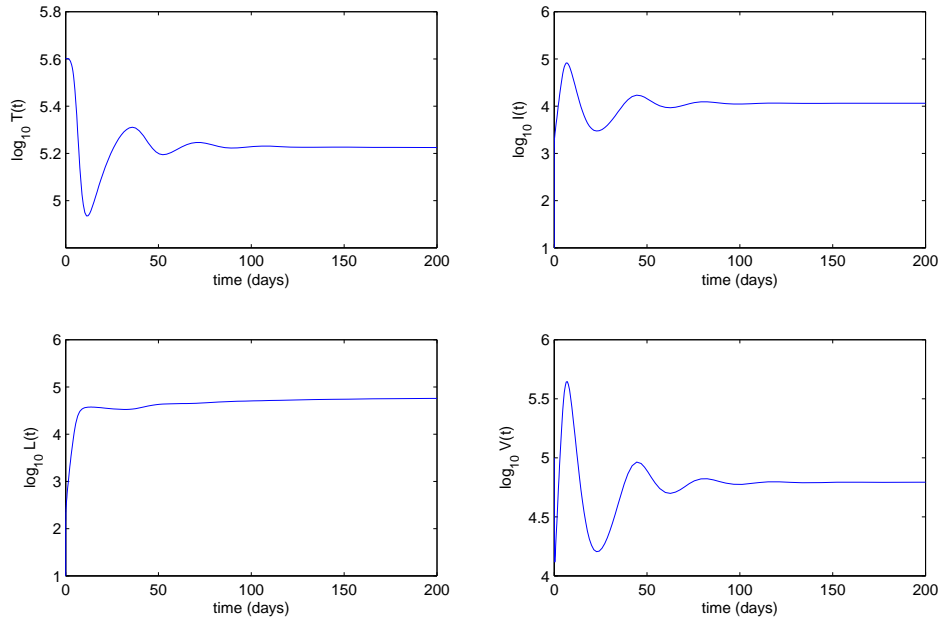}
\vspace{-0.2in}
\caption{A representative solution of (\ref{4CM}) with parameters from Table \ref{tab:1}.
The initial values are $T(0) = 4 \times 10^5$, $V(0) = 10^5$, and
$I(0) = L(0) = 0$.
In this example, the system tends to the endemic equilibrium as $t \to \infty$ because $R_L =2.027$. }
\label{fig:2}  
\end{center}
\end{figure*}
Figure \ref{fig:2} displays a representative graph of solutions to (\ref{4CM}) for which $R > 1$.  We note that the oscillations of $T, I$, and $V$ seem quite damped in comparison to those of Figure \ref{fig:1}.

\subsection{Model Parameters}

In this section and the previous one, we have adopted parameter values from other studies. 
A few of the parameters possess generally agreed upon values, including $\lambda$, $d_T$, $d_I$, and $d_V$.  However, it should be noted that $\lambda$ and $d_T$ are typically estimated for healthy individuals, and thus may not be reliable estimates to describe values within HIV patients, especially for those who experience impaired thymic function \cite{KimPerelson}. Obviously, there are many parameters, and these are summarized within Table 1, along with descriptions
of the variables, their associated units, and references from which parameter values stem.

The parameter that displays the most uncertainty within the literature is the viral infectivity $k$ which fluctuates by an order of magnitude from a value of $2.4 \times 10^{-8}$ ml/day \cite{RongPerelson, PKdB} to $2.7 \times 10^{-7}$ ml/day \cite{TuckLC}.  The value we utilize here is at the low end of this range and stems from \cite{PKdB}. 
Biologically relevant values of the \emph{in vivo} burst size $N$ are also somewhat uncertain.
Estimates based on counting HIV-1 RNA molecules in an infected cell vary between hundreds and thousands \cite{RongPerelson,HHZSF, HKDSS}, and estimates based on viral production have been as high as $5 \times 10^4$ \cite{CDPGH, dBRP}.
Here, we choose $N = 2000$ HIV-1 RNA/cell as reported in \cite{HKDSS}. 

Parameters that stem specifically from (\ref{4CM}) are generally not well-known. In particular, the fraction of new viral infections resulting in latency $\alpha$ varies from study to study, but based on previous work \cite{CallPerelson, KimPerelson}, we use $\alpha = 0.01$ per day. Similarly, the removal rate of latently infected cells, $d_L$, has been discussed as anywhere from $10^{-3}$ per day \cite{RongPerelson} to $0.24$ per day \cite{TuckLC}. Hence, we chose a value with this range, namely $d_L = 4 \times 10^{-3}$ per day, as reported within \cite{Finzi}.
The proportion $p$ of cells which are categorized as latent upon becoming infected also differs throughout recent studies, ranging from $1.5 \times 10^{-6}$ in \cite{CallPerelson} to $0.1$ in \cite{TuckLC}.  In this case, we utilize the latter value so as to emphasize changes in (\ref{3CM}) due to the latent cell population. 
Initial conditions of the proposed model are chosen to match clinically observed decay characteristics of the latent reservoir and the virus population. In particular, we utilize values similar to \cite{CallPerelson}, namely $T(0) = 4 \times 10^5$, $I(0)=0$, $L(0) = 0$, and $V(0) =  10^5$.
\begin{table*}[t]
\begin{center}
\normalsize
\caption{Variable and parameter values for (\ref{3CM}) and (\ref{4CM}).}
\label{tab:1}       
\begin{tabular}{lclcc}
\hline\noalign{\smallskip}
Variable & Units & Description & Value & Reference \\
\noalign{\smallskip}\hline\noalign{\smallskip}
$T(t)$ \qquad & cells ml$^{-1}$ \qquad & Susceptible CD$4+$ T-cells & -- & -- \\
$I(t)$ & cells ml$^{-1}$ & Actively Infected CD$4+$ T-cells & -- & -- \\
$L(t)$ & cells ml$^{-1}$ & Latently Infected CD$4+$ T-cells & -- & -- \\
$V(t)$ & virions ml$^{-1}$ & Infectious virions & -- & -- \\
\hline\noalign{\smallskip}
$\lambda$ &  ml$^{-1}$ day$^{-1}$ & Production rate of CD$4+$ T-cells \qquad & $10^4$ & \cite{CallPerelson}\\
$d_T$ & day$^{-1}$ & Death rate of susceptible $T$ cells & $0.01$ & \cite{Mohri} \\
$d_I$ & day$^{-1}$ & Death rate of actively infected $T$ cells & $1$ & \cite{Markowitz}\\
$d_V$ & day$^{-1}$ & Clearance rate of virions & $23$ & \cite{Ramratnam99}\\
$k$ & ml day$^{-1}$ & Rate of infection of susceptible cells & $2.4\times 10^{-8}$ & \cite{PKdB}\\
$N$ & - & Burst rate of actively infected T-cells & $2000$ & \cite{HKDSS} \\
\hline\noalign{\smallskip}
$d_L$ & day$^{-1}$ & Death rate of latent cells & $4 \times 10^{-3}$ & \cite{Finzi} \\
$\alpha$ & day$^{-1}$ & Activation rate of latent cells & $0.01$ & \cite{CallPerelson} \\
$p$ & - & Proportion of latent infection & $0.1$ & \cite{TuckLC}\\
\hline\noalign{\smallskip}
$\eps_{RT}$ & - & Efficacy of RT inhibitor & varies  & --   \\
$\eps_{PI}$ & - & Efficacy of protease inhibitor & varies   & -- \\
\noalign{\smallskip}\hline
\end{tabular}
\end{center}
\end{table*}
Next, we analyze properties of solutions to (\ref{4CM}) so as to compare their dynamics and large time behavior with solutions of (\ref{3CM}).

\subsection{Analysis and properties of solutions}

As a first step, we can say with certainty that biologically reasonable values of the parameters give rise to positive populations assuming that at some earlier point (perhaps at the initial time) the populations possessed positive values.
\begin{theorem}
\label{T1}
Assume all constants in (\ref{4CM}) are nonnegative and the initial values $T(0), I(0), L(0)$, and $V(0)$ are positive. Then, the solutions of (\ref{4CM}), namely $T(t),I(t),L(t)$, and $V(t)$ exist, are unique, and remain bounded on the interval $[0,t^*]$ for any $t^* > 0$.  Additionally, each function remains positive for any $t \geq 0$.
\end{theorem}
Of course, the requirement of initial positivity is not completely necessary since we may translate or rescale the time variable to alter the initial time.  Hence, what is necessary for the theorem to hold is that all populations must attain positive values at \emph{some} time.
This result provides some general validation for the model since it implies that negative population values cannot occur if one begins with biologically reasonable (i.e., positive) values.

Next, we proceed as for (\ref{3CM}) and investigate the possible equilibrium states of (\ref{4CM}) and their stability properties.
We find steady states by solving the nonlinear system of algebraic equations
\begin{equation}
\label{SS}
\left \{ \begin{aligned}
0 &= \lambda - d_T T - k T V \\
0 &=  (1-p)k T V + \alpha L - d_I I \\
0 &=  pk T V  - \alpha L - d_L L \\
0 &= Nd_I I - d_V V.
\end{aligned} \right.
\end{equation}
for the unknown constants $T, I, L$, and $V$ in terms of the parameters $\lambda, k, p$, $\alpha, N, d_T, d_I, d_L$, and $d_V$.
This is a nontrivial task, but eventually we find the existence of exactly two steady states.
We begin by solving for the nonlinear interaction term in the first equation so that
$$kTV = \lambda - d_T T.$$
With this, we may solve for $L$, $I$, and hence $V$ in terms of $T$ alone.  
From the third equation, $ pkTV = (d_L + \alpha) L$ and thus
$$ L = \frac{p}{d_L + \alpha} \left ( \lambda - d_T T \right).$$
Next, in the second equation, we find $(1-p) kTV = -\alpha L + d_I I$ and thus
$$I = \frac{1}{d_I} \left ( 1 -p + \frac{\alpha p}{d_L + \alpha} \right ) (\lambda - d_T T).$$
The last equation yields $V$ in terms of  $I$, whence $T$, so that
$$V = \frac{N d_I}{d_V} I = \frac{N}{d_V} \left ( 1 -p + \frac{\alpha p}{d_L + \alpha} \right ) (\lambda - d_T T).$$
Finally, we may use the representation of $V$ in terms of $T$ within the first equation and solve a simple quadratic in $T$
to determine the possible steady state values.  With this, the first equation becomes
$$0 = \lambda - d_T T - \frac{kN}{d_V} \left ( 1 -p + \frac{\alpha p}{d_L + \alpha} \right ) (\lambda - d_T T) T$$
and it follows that the only solutions are
$$ T = \frac{\lambda}{d_T} \qquad \mbox{and} \qquad T = \frac{\lambda}{d_T} \cdot \frac{1}{R_L}$$
where
\begin{equation}
\label{RL}
R_L = \frac{kN\lambda}{d_T d_V} \cdot \frac{(1-p) d_L + \alpha}{d_L + \alpha}.
\end{equation}
Continuing in this manner, we obtain two corresponding values for $I, L$, and $V$.
To summarize, we find two equilibria, which we write in the form $(T,I,L,V)$ as
$$\begin{gathered}
E_{NI}: \qquad  \left ( \frac{\lambda}{d_T}, 0, 0, 0 \right )\\
E_I: \qquad  \left ( \frac{\lambda}{d_T R_L},  \frac{d_T d_V}{kNd_I}(R_L - 1), \frac{p\lambda}{R_L(d_L + \alpha)} (R_L - 1), \frac{d_T}{k} (R_L - 1) \right ).
\end{gathered}$$

As before, we denote the non-infective equilibrium by $E_{NI}$ and the infective equilibrium by $E_I$. Notice that the limiting values of $T, I$, and $V$ for the infective state are of the same form as those of (\ref{3CM}), with $R_L$ replacing the role of $R_0$. Additionally, we see that if $R_L = 1$, then the equilibria coincide, and if $R_L < 1$, then the endemic equilibrium corresponds to negative values which, in view of Theorem \ref{T1}, cannot be obtained from biologically relevant initial data.

By studying the linearized version of the system, we may examine the local stability of these equilibria and find that their behavior mimics that of (\ref{3CM}).

\begin{theorem}
\label{T2}
If $R_L \leq 1$, then the non-infective equilibrium is locally asymptotically stable.  If $R_L > 1$ then the non-infective equilibrium is an unstable saddle point, and the endemic equilibrium is locally asymptotically stable.
\end{theorem}
Therefore, if $R_L \leq 1$ and population values begin within a sufficiently close distance of $E_{NI}$, then they will tend to $E_{NI}$ as $t \to \infty$. Contrastingly, if $R_L > 1$ and initial populations are sufficiently close to $E_I$, they will tend to $E_I$ in the long run.
Theorem \ref{T2} also emphasizes the crucial feature that equilibria are not stable simultaneously, that is, bistability of $E_{NI}$ and $E_I$
does not occur. Furthermore, it expresses that the qualitative behavior of system (\ref{4CM}) changes exactly when $R_L$ transitions from less than one to greater than one, and hence a bifurcation occurs at $R_L = 1$.

The final theorem of the section demonstrates the stronger result that initial values of these populations have no effect on their long term ($t \to \infty$) limiting values.
\begin{theorem}
\label{T3}
If $R_L \leq1$, then the non-infective equilibrium is globally asymptotically stable.  If $R_L > 1$, then the endemic equilibrium is globally asymptotically stable.
\end{theorem}

This analysis reveals one very important fact about the overall system: the end states of populations are only dependent on the value of $R_L$, and not any other parameter or initial value. If $R_L$ is greater than one, then the system tends to $E_{I}$, an end state with a non-zero population of infected cells and virions, but if $R_L$ is less than one, then the final equilibrium is $E_{NI}$, which contains neither virions nor infected T-cells.

The most important feature of these results is the explicit formula for $R_L$, which can be related exactly to the basic reproductive number of the three-component model (\ref{3CM}).  In order to investigate the differences between the two reproductive ratios, we define the quantity
\begin{equation}
\label{Q}
Q := \frac{R_L}{R_0} = \frac{(1-p) d_L + \alpha}{d_L + \alpha}.
\end{equation}
Notice that $Q$ depends only upon the three new parameters included within (\ref{4CM}), namely the activation ratio $\alpha$, proportion of cells which become latent upon infection $p$, and the death rate of latent cells $d_L$. Additionally, if the proportion $p$ of infections which produce latently infected T-cells is identically zero, then $R_L = R_0$.  However, since we consider  $p \in (0,1)$ we find
$$Q < \frac{d_L + \alpha}{d_L + \alpha} = 1$$ and the relationship $$R_L  < R_0$$ follows directly.
Thus, the reproduction number of the latent cell model (\ref{4CM}) is \emph{strictly less} than that of the standard three-component model (\ref{3CM}).  
Therefore, the stability of the the non-infective state is enhanced by the inclusion of the latently-infected cell population.  
Namely, there are more values of $\lambda, d_T, d_I, d_V, k$ and $N$ which correspond to $R_L \leq 1$ than $R_0 \leq 1$.
From a modeling standpoint, this result is somewhat intuitive.  Because (\ref{4CM}) assumes that a fraction of newly infected cells become latently infected and the latter can only activate (becoming actively productive) or die, the average number of infected cells generated by the introduction of a single infected cell into a susceptible system is decreased in comparison to a model without latently infected cells, namely (\ref{3CM}).  Hence, one should expect that the basic reproduction number, representing this average number of infected cells, does in fact decrease. 
Another consequence of this results is that there exist a number of parameter values for which $R_0 > 1$ but $R_L \leq 1$,
and in such cases the solutions of (\ref{3CM}) tend to $E_{NI}$ while those of (\ref{4CM}) tend to $E_I$ as $t \to \infty$.  
Clearly, the converse ($R_0 \leq 1$ but $R_L > 1$) is not possible by the above inequality.  
 In fact, we may rewrite their ratio $Q$ as
$$ Q = 1 - \frac{p d_L}{d_L + \alpha}$$ so that the difference between $R_0$ and $R_L$ is greatest for large values of $p$ and $d_L$, but small values of $\alpha$.  With the representative parameter values given in Table \ref{tab:1}, we see that 
$$R_L \approx 5.978 \quad \mbox{and} \quad R_0 \approx 6.154.$$ 
Hence, the change in system behavior caused by the difference between the reproduction numbers appears somewhat negligible, as both values are significantly larger than their respective bifurcation points.
Additionally, $Q = 0.97$ in this case, so that the relative difference between $R_0$ and $R_L$ is merely
$$\frac{R_0 - R_L}{R_0} = 1 - Q = 3\%.$$
Exactly quantifying this relative change, however, is difficult since many of the parameter values of Table \ref{tab:1},
in particular $k, p, \alpha$, and $d_L$, are not well-established, and hence this percentage could be much larger or perhaps even smaller.
For example, if we utilize the smallest value of $\alpha = 3 \times 10^{-3}$, stemming from \cite{PKdB}, and the largest  values of $p = 0.1$ \cite{TuckLC} and $d_L = 0.24$ \cite{PKdB}, then a simple computation shows that $Q = 0.9$. 
Hence, the relative difference between reproductive ratios could possibly be as large as $10\%$, 
though more conservative estimates of the associated parameters place the value of $1-Q$ between $1\%$ and $5\%$.

Regardless of the quantified distinction between the reproductive ratios, it seems somewhat counterintuitive that $R_L < R_0$, especially since so many authors \cite{Arlen, Blankson, CF, Chun1995, Chun1997a, Chun1997b, KimPerelson, KPAIDS, PECVH, RongPerelson, SmithAgg} have detailed the large degree to which latent reservoirs contribute to the increased persistence of HIV infection via viral rebound after treatment with ART.  Hence, the result of the mathematical analysis, namely that the effects of latent infection \emph{reduce} the basic reproductive ratio, seems to contradict this theory.  However, as we previously stated, the alterations in the mathematical model explain this effect.  Additionally, the reproductive ratio is but one parameter, and so it seems unlikely that this particular metric will completely determine the realistic behavior of the system.  In fact, a more detailed analysis of the behavior of solutions over the timescales of biological relevance, rather than considering only the limiting behavior as $t \to \infty$, will demonstrate the shortcomings of the basic reproductive ratio.
We illustrate this using the effects of antiretroviral therapy and some associated computational results within the next section.

\section{Antiretroviral Therapy}

In order to further contrast these two models and the effects of the latent cell population, we will introduce additional parameters
to represent the application of antiretroviral therapy.  The inclusion of ART will allow us to determine the range of drug efficacies that
distinguish between the limiting dynamics of (\ref{3CM}) and (\ref{4CM}) and further elucidate the differing behaviors of the two models.

Two classes of antiretroviral drugs are often used to reduce the viral load and
limit the infected T-cell population. 
One class is known as Reverse Transcriptase Inhibitors (RTIs), which can block
the infection of target T-cells by infectious virions. 
The other category is Protease Inhibitors (PIs), which prevent HIV-1 protease from cleaving the HIV polyprotein into
functional units, thereby causing infected cells to produce immature virus particles
that are non-infectious. 
In this way, RTIs serve to reduce the rate of infection of activated CD$4+$ T-cells, whereas PIs
decrease the number of new infectious virions that are produced. Both drugs thus diminish the propagation
of the virus \cite{Janeway, RongPerelsonRev2}. 
While we expect that latently infected cells may absorb PIs and that such
cells, when activated, will produce noninfectious virus, we 
will instead assume that PIs have no effect on the proportion of cells
that are latently infected. 
This is in line with some experimental findings, that suggest
that antiretroviral drugs do not effectively block replication of virus from the latent reservoir \cite{Chun2003}. 
Hence, in our model, susceptible T-cells may be inhibited with either RTIs, or PIs, or they may become infected. Infected
cells may be inhibited with PIs, and cells inhibited with one drug may be inhibited with the other. 
In the presence of these two inhibitors, the model equations (\ref{4CM}) are modified to become:
\begin{equation}
\label{4CMART}
\left \{ \begin{aligned}
\frac{dT} {dt} &= \lambda - d_T T - k(1- \eps_{RT}) T V_I \\
\frac{dI} {dt} &=  (1-p)k(1- \eps_{RT}) T V_I + \alpha L - d_I I \\
\frac{dL} {dt} &=  pk(1- \eps_{RT}) T V_I  - \alpha L - d_L I \\
\frac {dV_I} {dt} &= N(1- \eps_{PI})d_I I - d_V V_I.
\end{aligned} \right.
\end{equation}
where $\eps_{RT}, \eps_{PI} \in [0,1]$ are the efficacies of RTIs and PIs, and $V_I$ represents 
the population of infectious virions.  We may include the number of non-infectious
virions $V_{NI}$, with the total viral load $V = V_I + V_{NI}$, but $V_{NI}$
decouples from the remaining equations, and hence plays no role in the evolution of the system.

\begin{figure*}[t]
\begin{center}
\includegraphics[scale = 0.6]{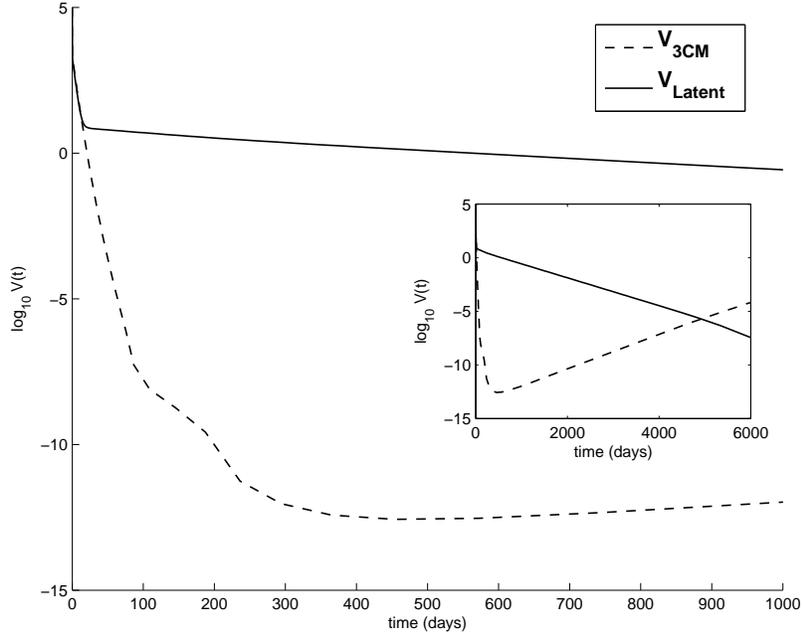}
\vspace{-0.2in}
\caption{Viral loads for (\ref{4CMART}) and (\ref{3CMART}) with $\eps = 0.519$.
This value of $\eps$ yields $R^\eps_0 = 1.003$ and $R^\eps_L = 0.974$, respectively.
The inset figure includes the viral behavior for $t \in [0,6000]$ days.}
\label{fig:3}  
\end{center}
\end{figure*}

To compute the steady states and basic reproduction number for (\ref{4CMART}), we may reproduce the analysis of (\ref{4CM}), but clearly the new terms are introduced only where the parameters $k$ and $N$ appear.  Thus, we need only replace $k$ with $k(1- \eps_{RT})$ and $N$ with $N(1- \eps_{PI})$.
The new basic reproduction number then becomes 
\begin{equation}
\label{RLART}
R^\eps_L = \frac{kN(1-\eps)\lambda}{d_T d_V} \cdot \frac{(1-p) d_L + \alpha}{d_L + \alpha}
\end{equation}
where we define the quantity $\eps = \eps_{RT} + \eps_{PI} - \eps_{RT}\eps_{PI}$, so that
$$1- \eps = (1-\eps_{RT})(1- \eps_{PI}).$$
Since (\ref{4CMART}) is identical to (\ref{4CM}) with the minor change in parameter values described above, Theorems \ref{T2} and \ref{T3} hold for (\ref{4CMART}) with the corresponding value $R^\eps_L$ instead of $R_L$.
Writing the corresponding viral steady state for (\ref{4CMART}) we find $$\overline{V} = \frac{d_T R_L}{k}(1 -\eps_{PI}) - \frac{d_T}{k(1- \eps_{RT})}$$ and we notice that its partial derivative
$$\frac{\partial \overline{V}}{\partial \eps_{RT}} = - \frac{d_T}{k(1-\eps_{RT})^2}$$ is large, especially when $\eps_{RT} \approx 1$. Thus,
$\overline{V}$ is sensitive to small changes in $\eps_{RT}$, and this sensitivity increases with the efficacy of the RTI.
Hence, this model does not realistically describe the persistence of low-level viremia in patients on reverse transcriptase inhibitors, as previously
addressed within \cite{RongPerelson}.  However, we note that 
$$\frac{\partial \overline{V}}{\partial \eps_{PI}} = - \frac{d_T R_L}{k}$$ which is constant for all values of $\eps_{PI}$, and does not possess the same sensitivity.
Thus, to simplify the analysis, we will assume throughout that only PIs are used, and hence $\eps = \eps_{PI}$ while $\eps_{RT} = 0$.

Upon incorporating the use of antiretroviral drugs into (\ref{3CM}), the system becomes
\begin{equation}
\label{3CMART}
\left \{ \begin{aligned}
\frac{dT} {dt} &= \lambda - d_T T - k(1- \eps_{RT}) T V_I \\
\frac{dI} {dt} &=  k(1- \eps_{RT}) T V_I - d_I I \\
\frac {dV_I} {dt} &= N(1- \eps_{PI})d_I I - d_V V_I.
\end{aligned} \right.
\end{equation}
with associated basic reproduction number 
\begin{equation}
\label{R0ART}
R_0^\eps = \frac{kN(1-\eps)\lambda}{d_T d_V}.
\end{equation}
As for (\ref{4CMART}), the stability results for (\ref{3CM}) contained in Section $2$ hold for (\ref{3CMART}) by replacing $R_0$ with $R^\eps_0$.
Comparing the two values $R^\eps_0$ and $R^\eps_L$, we see that their ratio is again $Q$ given by (\ref{Q}), so that $R^\eps_L = Q R^\eps_0$.  
As before, since $Q < 1$ we find $R^\eps_L < R^\eps_0$ and even with the incorporation of ART, 
the latent cell population \emph{decreases} the basic reproduction number of the system.  
With the representative parameter values given in Table \ref{tab:1}, we see that 
the basic reproduction numbers associated with (\ref{4CMART}) and (\ref{3CMART}) are
$$R^\eps_L \approx 2.027(1-\eps) \quad  \mbox{and} \quad R^\eps_0 \approx 2.087(1-\eps).$$
Hence, in order for the non-infective state to be realized in (\ref{3CMART}), we must
have $\eps > 0.521$, which only differs mildly from the value of $\eps$ that is
needed to reach the same non-infective state in (\ref{4CMART}), namely $\eps >  0.506$.
Thus, the antiretroviral therapy must attain an efficacy only $3\%$ greater in order to achieve the analogous effect.
From this analysis, it would seem that the establishment of a latent reservoir should not strengthen a continued
infection since the drug efficacy necessary to drive the system to a non-infective equilibrium is actually less 
for (\ref{4CMART}) than for (\ref{3CMART}).  
This result seems to greatly contradict the known issues that scientists have faced regarding the eradication
of the viral reservoir.  
However, the resolution of these seemingly opposing viewpoints is made quite clear by precise numerical simulations.

As an illustrative example of the difference in long-time dynamics between the two models, we may choose values of $\eps$ which yield  $R^\eps_0 >1$ and $R^\eps_L < 1$ and measure their corresponding behavior.  A representative simulation is presented in Figure \ref{fig:3}.  To differentiate between the corresponding viral loads, we will denote the infectious virus population associated with (\ref{3CMART}) by $V_{3CM}(t)$ and its latent model analogue by $V_{Latent}(t)$ as in the figure.  In this case, we choose $\eps = 0.519$ and find $R^\eps_0 = 1.003$ and $R^\eps_L = 0.974$.  Hence, from the known results of Section $2$ and Theorem (\ref{T3}), we may deduce that $V_{3CM}(t) \to \frac{d_T}{k} (R^\eps_0 - 1)  \approx 2 \times 10^3$ as $t \to \infty$, while $V_{Latent}(t) \to 0$ as $t \to \infty$. However, one can distinctly see from Figure \ref{fig:3}
that the early decay rate of $V_{3CM}$ is much greater than that of $V_{Latent}$.

Notice that the effects of latent infection do not influence the viral load for the first thirty days of treatment as $V_{3CM}(t)$ and $V_{Latent}(t)$ follow the same approximate trajectory during this time period.  However, once the latently infected T-cell population grows sufficiently large, the effects are tremendous.  Throughout the first three years of continuous treatment, $V_{3CM}$ diminishes greatly, past $10^{-10}$ in fact, while $V_{Latent}$ remains $O(1)$ even up to day $1000$.  Certainly this seems strange as $R^\eps_0 > 1$ implies a persistent virus population must develop for (\ref{3CMART}) and $R_L < 1$ dictates the eventual elimination of the viral population for the latent model.  Within the inset figure, it can be seen that the behavior predicted by the basic reproduction numbers does eventually occur, that is, values of $V_{3CM}$ rebound and tend to a positive equilibrium, while those of $V_{Latent}$ continue their slow, steady decline to eradication.  Unfortunately, these events occur nearly fifteen years after the introduction of ART and well outside the timescale of biological relevance.  Hence, it appears that the values of $R^\eps_0$ and $R^\eps_L$ alone cannot provide sufficient information to account for the realistic biological dynamics of the model due to the change in timescales and decay rates introduced by the latent cell population.  A better estimate of the behavior would certainly be given by the rates of decay to eradication, but precise estimates on these quantities are more difficult to obtain analytically.  Instead, we examine a slightly different metric of viral persistence or clearance. 

\begin{figure*}[t]
\begin{center}
\includegraphics[scale = 0.58]{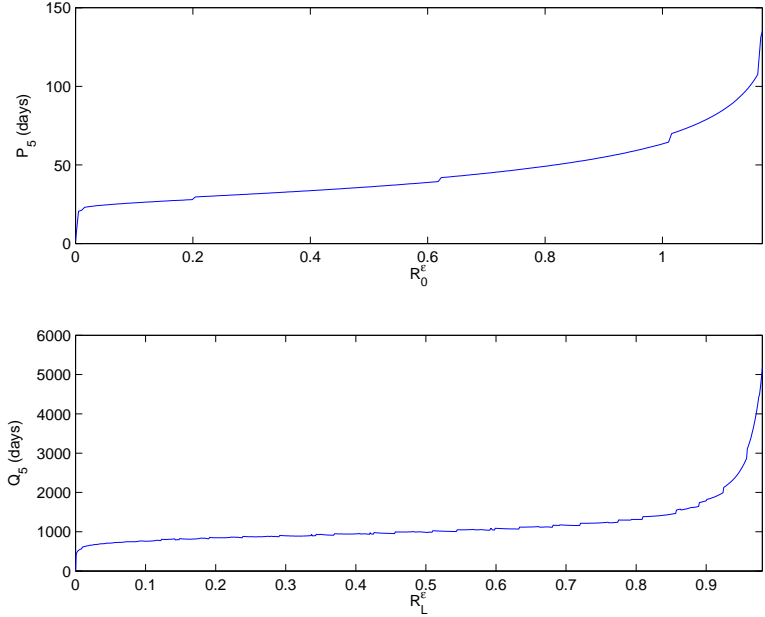}
\vspace{-0.2in}
\caption{Comparison of $P_5(R_0^\eps)$ and $Q_5(R_L^\eps)$ }
\label{fig:4}  
\end{center}
\end{figure*}

The feature of viral clearance that one must capture here is not just the decay of the viral load, but a sufficiently rapid speed of decay so as to be realized within a time period of biological relevance.  Hence, we consider a specific value of the virus population to represent clearance, and proceed to study the minimum arrival time of the viral load to that value.
In this vein, we define the functions $$P_n(r) = \inf \left \{t > 0 : \log_{10} \biggl (V_{3CM}(t) \biggr) \leq -n \ \mbox{for} \ R^\eps_0 = r \right \}$$ and
$$Q_n(r) = \inf \left \{t > 0 : \log_{10}\biggl (V_{Latent}(t) \biggr ) \leq -n  \ \mbox{for} \ R^\eps_L = r \right \}.$$
We note that either of these functions may become infinite if the population of virions fails to reach the prescribed value for any positive time.  For example, $P_5(2) = Q_5(2) = \infty$, since neither viral load obtains a value as small as $10^{-5}$ for a corresponding reproductive ratio of $2$, while $P_5(0.6) \approx 25$ and $Q_5(0.6) \approx 1000$ as represented in Figure \ref{fig:4}.  Unlike the stability of equilibria, the values of $P_n(r)$ and $Q_n(r)$ will depend upon the initial population values that are chosen.  Within the present study, however, we will continue to utilize the initial populations of previous sections to serve as a representative example.

We first select the value of $10^{-5}$ copies per ml for our definition of viral eradication and study the associated times to eradication provided by the functions $P_5$ and $Q_5$. 
Namely, what we are assuming is that once the viral population is suitably dilute - in this case less that $10^{-5}$ copies per ml - then the infection has been cleared and no rebound can occur.  Figure \ref{fig:4} provides a comparison of $P_5$ and $Q_5$ for differing values of $R^\eps_0$ and $R^\eps_L$, respectively.  Though their general shapes are quite similar, the associated time periods differ dramatically.  Typical values of $P_5$ range from $20$ to $100$ days, while the majority of values of $Q_5$ range between $1000$ and $3000$ days.  As can be seen in Figure \ref{fig:4}, even if the efficacy of the RT inhibitor, $\eps$, approaches $100\%$, and thus $R^\eps_0$ approaches zero, $V_{3CM}$ requires around $15$ days to reach a value of $10^{-5}$.  In this same situation, $R^\eps_L$ approaches zero, but $V_{Latent}$ requires nearly $500$ to $1000$ days to reach a value of $10^{-5}$.  Thus, even for values of $R^\eps_L$ which are significantly less than one, we see that it would require nearly three years for the viral load to reach this threshold due to the influence of latent infection.  In addition, we see that $V_{3CM}$ will reach values of $10^{-5}$ even if $R^\eps_0 > 1$, and this will occur within $100$ days, almost ten times faster than it would take $V_{Latent}$ to reach the same value for a constant drug efficacy around $90\%$.

\begin{figure*}[t]
\begin{center}
\includegraphics[scale = 0.58]{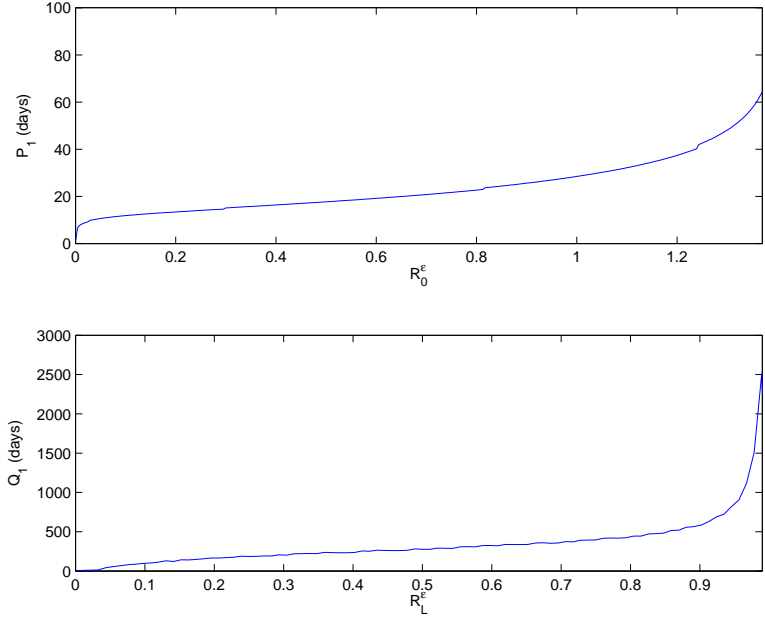}
\vspace{-0.2in}
\caption{Comparison of $P_1(R_0^\eps)$ and $Q_1(R_L^\eps)$ }
\label{fig:5}
\end{center}
\end{figure*}

Considering that the biological detection threshold is around $50$ viral copies per ml \cite{Doyle, CallPerelson}, one possibility is that the $10^{-5}$ threshold above has been chosen too small in Figure \ref{fig:4} to effectively serve as a realistic measure of eradication.  Hence, we consider $P_1$ and $Q_1$ and perform a similar analysis.  Figure \ref{fig:5} contains these simulations and displays a decrease in the time necessary to reach the defined threshold of $10^{-1}$.  However,  even for a $70\%$ constant drug efficacy, in which case $R^\eps_L = 0.6$, we see that approximately one year of continuous ART would still be required to reach a viral load of $0.1$ copies per ml for (\ref{4CMART}).
Additionally, values of $P_1$ remain around $60$ for (\ref{3CMART}) even if $R^\eps_0$ is near $1.5$, which exceeds the bifurcation point by nearly $50\%$.  Thus, if we define viral clearance as a decay in the viral load to $0.1$ copies per ml within six months of treatment, then (\ref{3CMART}) would require $R^\eps_0$ to be less than $1.5$ while (\ref{4CMART}) would require $R^\eps_L$ to be less than $0.2$. Obviously, a much wider range of parameter values will yield $R^\eps_0 < 1.5$ than $R^\eps_L < 0.2$, and we see that the latent T-cell population does, in fact, extend the period of time during which viremia persists, even though the behavior as $t \to \infty$, as given by Theorems \ref{T2} and \ref{T3}, may provide seemingly contrary information.

From this, the biological influence of latent infection becomes clear - the time needed to decrease the viral load to values from which rebound is unlikely or unable to occur is increased by a factor of ten or twenty. Hence, the value of the basic reproduction number alone does not represent a proper definition for viral persistence or eradication, and the functions provided above $P_n(r)$ and $Q_n(r)$, for well-chosen values of $n$, possess the information required to better determine the behavior of the infection. Further analysis can be performed for smaller (and negative) values of $n$, but the results discussed above are typical. 
In the next section, we prove Theorems \ref{T1}, \ref{T2}, and \ref{T3} regarding the qualitative behavior of the latent infection model.

\section{Proofs of main theorems}
With the analysis concluded, we finally prove the main results of the previous sections. In what follows, $C$ will be used to denote a positive, but arbitrary constant which may change from line to line. First, we prove the existence, uniqueness, and positivity of solutions.
\begin{pf}[Theorem \ref{T1}]
While one may prove that a certain positive set remains invariant under the flow (as in \cite{Elaiw}), this requires assumptions which bound the initial data from above.  In our proof, we utilize a continuity argument instead and do not assume any upper bounds on initial data.
Using the Picard-Lindeloff theorem and the quadratic nature of the equation, the local-in-time existence of a unique, $C^1$ solution follows immediately.
Hence, we will concentrate on proving positivity of solutions as long as they remain continuous, and this property will yield bounds on the growth of solutions. 
From the bounds obtained below, then, it follows that the solution exists globally and is both unique and continuously differentiable for all $t > 0$.  
Define $$T^* = \sup \{t \geq 0 : T(s), I(s), L(s), V(s) > 0, \ \mbox{for all} \ s \in [0,t] \}.$$ 
Since each initial condition is nonnegative and the solution is continuous, there must be an interval on which the solution remains positive, and we see that $T^* > 0$.
Then on the interval $[0,T^*]$ we estimate each term.

Lower bounds on $I, L$, and $V$ instantly follow since the decay terms are linear. 
More specifically, we find
$$\frac{dI}{dt}=(1-p)kTV + \alpha L - d_I I \geq -d_I I$$ and thus
$$ I(t) \geq I(0) e^{-d_I t} > 0$$ for all $t \in [0,T^*]$.
Similarly, for the latent T-cell population
$$\frac{dL}{dt}=pkTV - (d_L + \alpha) L \geq - (d_L + \alpha) L$$ and thus
$$ L(t) \geq L(0) e^{-(d_L + \alpha) t} > 0$$ for all $t \in [0,T^*]$.
The positivity of the virion population follows in the same manner since
$$\frac{dV}{dt}= N d_I I - d_V V \geq - d_V V$$ and thus
$$ V(t) \geq L(0) e^{-d_V t} > 0$$ for all $t \in [0,T^*]$.
The positivity of $T$ requires extra effort since it decreases due to the nonlinearity.
We first construct an upper bound on $\frac{dT}{dt}$ as
$$\frac{dT}{dt}= \lambda-\mu T -kTV \leq \lambda$$
and thus
$$ T(t) \leq T(0) + \lambda t \leq C(1+t).$$
Next, we sum the equations for $I, L$, and $V$, and by positivity of these functions, obtain upper bounds on each one. Using the upper bound on $T(t)$, we find
$$\frac{d}{dt}(I + L + V)=kTV  + (N-1)d_I I - d_L L - d_V V \leq C(1+t) \left ( I + L + V \right ). $$
By Gronwall's Inequality, we have
$$I(t) + L(t) + V(t)  \leq Ce^{t^2}$$
for $t \in [0,T^*]$.
Since $I(t)$ and $L(t)$ are positive on this interval, the same upper bound follows on $V(t)$ alone.  
With this, we can now obtain a lower bound on $T$.  We find
$$\frac{dT}{dt}=\lambda- d_T T -kTV \geq - d_T T - kTV \geq -C(1+e^{t^2})T$$
or stated equivalently
$$\frac{dT}{dt} + C(1+e^{t^2})T \geq 0.$$
It follows that 
$$\frac{d}{dt} \left ( T(t)e^{C\int_0^t(1+e^{\tau^2})d\tau} \right ) \geq 0$$
and $T(t) \geq T(0) e^{-C\int_0^t(1+e^{\tau^2})d\tau} > 0$
for $t \in [0,T^*]$.
Finally, if $T^* < \infty$, then all functions are strictly positive at time $T^*$, contradicting its definition as the supremum of such values.  Hence, we find $T^* = \infty$ and the result follows.
\end{pf}

Next, we prove the local stability results.
\begin{pf}[Theorem \ref{T2}]

We proceed by linearizing the system and using the Routh-Hurwitz criterion to determine conditions under which the linear system possesses only negative eigenvalues. Then, as a consequence of the Hartman-Grobman Theorem, the local behavior of the linearized system is equivalent to that of the nonlinear system.

First, we compute the Jacobian evaluated at the non-infective equilibrium $E_{NI} = (\frac {\lambda} {d_T} , 0 ,0,  0)$, resulting in
$$J(E_{NI})  = \begin{bmatrix}
	\ - d_T \ & \ 0 \ & \ 0 \ & \ - \frac {k \lambda} {d_T} \ \\[4 pt]
	0 & -d_I &  \alpha & \frac {k(1-p) \lambda} {d_T}  \\[4 pt]
	0 & 0 & - d_L - \alpha & \frac {kp \lambda} {d_T}  \\[4 pt]
	0 & N d_I & 0 & -d_V
\end{bmatrix}$$
From this, we compute the associated characteristic polynomial for eigenvalues $\eta$
\begin{eqnarray*}
0 & = & \det( J - \eta \mathbb{I}) \\
& = &  (d_T + \eta) \left [(d_I + \eta)(d_L + \alpha + \eta)(d_V + \eta) - \frac{\alpha N K \lambda p d_I}{d_T} \right.\\
& \ & \quad  \left. - \ \frac{N K \lambda (1-p) d_I}{d_T} \left (d_L + \frac{\alpha}{1-p} + \eta \right ) \right ].
\end{eqnarray*}
After expanding the terms and ordering by powers of $\eta$, this equation ultimately simplifies to
\begin{equation}
\eta^3+A_1\eta^2+A_2\eta+A_3  = 0
\end{equation}
where 
$$\begin{gathered}
A_1 =  d_V + d_I + d_L + \alpha\\
A_2 =  d_I d_V + (d_L + \alpha)(d_I + d_V) - \frac{(1-p) \lambda N k d_I}{d_T}  \\
A_3 = (d_L + \alpha)d_I d_V - \frac{\lambda N k d_I}{d_T} ((1-p) d_L + \alpha).
\end{gathered}$$
The Routh-Hurwitz criterion requires $A_1, A_2, A_3 > 0$ and  $A_1A_2 - A_3 > 0$.  Clearly, $A_1 > 0$, and after rewriting $A_3$ in terms of $R_L$, we find $$A_3 = (d_L + \alpha)d_I d_V (1- R_L).$$ 
Thus, in order for all of the eigenvalues of the system to be negative, it is necessary that $R_L < 1$.  Similarly, we rewrite $A_2$ as $$A_2 = (d_L + \alpha)(d_I + d_V) + d_I d_V \left [ 1 - R_L\frac{(1-p) (d_L + \alpha)}{(1-p) d_L + \alpha} \right ].$$  
Using the inequality
\begin{equation}
\label{Ineq}
\frac{(1-p) (d_L + \alpha)}{(1-p) d_L + \alpha} = 1 - \frac{p \alpha}{(1-p) d_L + \alpha} < 1,
\end{equation}
and the previous condition $R_L < 1$, we find $A_2 > 0$.

Finally, using (\ref{Ineq}), we see that $A_2 > d_I d_V(1 - R_L)$, and clearly $A_1 > d_L + \alpha$.  Therefore, we find
$$A_1A_2 >  d_I d_V(d_L + \alpha)(1- R_L) = A_3$$
and the Routh-Hurwitz criteria are satisfied.
Thus, $R_L< 1$ implies that all eigenvalues of the linearized system are negative, and hence the local asymptotic stability of $E_{NI}$ follows.
Conversely, if $R_L >1$, then the linearized system possesses at least one positive eigenvalue, and the equilibrium is unstable.\\

The analysis for $E_I$ is similar to that of $E_{NI}$. For notational purposes, we label the equilibrium as $(\bar{T}, \bar{I}, \bar{L}, \bar{V})$. Linearizing (\ref{4CM}) about $E_I$, we find the Jacobian 
$$J(E_I)  = \begin{bmatrix}
	\ - (d_T + k \bar{V}) \ & \ 0 \ & \ 0 \ & \ - k \bar{T} \ \\[4 pt]
	(1-p) k \bar{V} & -d_I &  \alpha & k(1-p) \bar{T}  \\[4 pt]
	pk \bar{V} & 0 & - d_L - \alpha & kp \bar{T}  \\[4 pt]
	0 & N d_I & 0 & -d_V \\
\end{bmatrix}$$
and this results in the characteristic equation
\begin{eqnarray*}
0 & = & (d_T R_L + \eta) (d_I + \eta)(d_L + \alpha + \eta)(d_V + \eta)\\
& \ & \qquad - \frac{N K \lambda(1-p) d_I}{d_T R_L} \left ( d_L + \frac{\alpha}{1-p} + \eta \right )(d_T + \eta).
\end{eqnarray*}
After expanding terms and simplifying, we arrive at the quartic polynomial
$$\eta^4 + A_1\eta^3+A_2\eta^2+A_3\eta+A_4  = 0$$
where 
\begin{align*}
A_1 & =  d_T R_L + d_V + d_I + d_L + \alpha\\
A_2 & =  d_TR_L (d_L + \alpha + d_I + d_V) + (d_L + \alpha)(d_I + d_V) + d_Id_V \\
 & \qquad - \frac{(1-p) \lambda N k d_I}{d_T R_L}  \\
A_3 & = d_TR_L(d_L + \alpha)(d_I + d_V) + d_TR_Ld_I d_V + (d_L + \alpha)d_Id_V \\
& \qquad - \frac{\lambda N k d_I}{d_TR_L} ((1-p)d_T +  (1-p)d_L + \alpha)\\
A_4 & = d_TR_L(d_L + \alpha)d_I d_V - \frac{\lambda N k d_I}{R_L} ((1-p) d_L + \alpha)
\end{align*}
As before, the Routh-Hurwitz criterion requires all coefficients to be positive, as well as,  $A_1A_2 - A_3 > 0$ and now $A_3(A_1A_2 - A_3) - A_4A_1^2 > 0$.  As for the $E_{NI}$ analysis, the positivity of $A_1$ follows directly from the positivity of the coefficients, and after rewriting $A_4$, we find $$A_4 = d_T(d_L + \alpha)d_I d_V (R_L-1).$$ 
Hence, it is necessary that $R_L > 1$ in order to satisfy the criteria.  Similarly, we rewrite $A_3$ as
\begin{eqnarray*}
A_3 & = & d_TR_L(d_L + \alpha)(d_I + d_V) + d_TR_Ld_I d_V + (d_L + \alpha)d_Id_V\\
& \ & \qquad - \left [ d_Td_Id_V \frac{(1-p)(d_L + \alpha)}{(1-p) d_L + \alpha} + d_Id_V(d_L + \alpha) \right ] \\
& > & d_TR_L(d_L + \alpha)(d_I + d_V) + d_Td_Id_V(R_L - 1) > 0. 
\end{eqnarray*}  
In this inequality we have canceled the third term with the last term and utilized the inequality (\ref{Ineq}) to bound the fourth term.
The only nonpositive term in $A_2$ can be rewritten as $$ - \frac{(1-p) \lambda N k d_I}{d_T R_L}  = -d_Id_V \frac{(1-p)(d_L + \alpha)}{(1-p) d_L + \alpha} > -d_Id_V,$$ and hence $$A_2 > d_TR_L (d_L + \alpha + d_I + d_V) + (d_L + \alpha)(d_I + d_V) > 0.$$

By the definition of $A_1$, we have $A_1 > d_I + d_V$ and using the above inequality for $A_2$, we find
\begin{eqnarray*}
A_1 A_2 & > & (d_I + d_V) \cdot \left [d_TR_L (d_L + \alpha + d_I + d_V) + (d_L + \alpha)(d_I + d_V) \right ]\\
& > & (d_I + d_V) \cdot d_T R_L (d_L + \alpha) + d_V \cdot d_T R_L d_I + d_V \cdot (d_L + \alpha) d_I \\
& > & A_3.
\end{eqnarray*}

Finally, we verify the last inequality, namely $A_3(A_1A_2 - A_3) - A_4A_1^2 > 0$.  After a long calculation, we find
\begin{eqnarray*}
A_1A_2 - A_3 & > & (d_L + \alpha) \biggl [ d_T R_L (d_T R_L + 2d_I + 2d_V + d_L + \alpha) \\
& \ & + (d_L + \alpha)(d_I + 2d_V) +  (d_I + d_V)^2    \biggr ].
\end{eqnarray*}
Removing selected terms from the right side, we also find
\begin{equation}
\label{123}
A_1A_2 - A_3 > d_I \biggl [ (d_L + \alpha)^2 + (d_L + \alpha)d_I + d_TR_L(d_L + \alpha) \biggr ].
\end{equation}
In addition, we see from a previous computation that
$$A_3 > d_Td_VR_L(d_L + \alpha) + d_Td_Id_V(R_L - 1).$$
Hence, we obtain a lower bound for $A_3(A_1A_2 - A_3)$ by multiplying the first term in the inequality for $A_3$ by the
right side of (\ref{123}) and the second term of the $A_3$ inequality by the previous lower bound for $A_1A_2 - A_3$.
This results in 
\begin{eqnarray*}
A_3 \left (A_1A_2 - A_3 \right ) & > & d_Td_VR_L(d_L + \alpha) \cdot d_I \biggl [ (d_L + \alpha)^2 + (d_L + \alpha)d_I + d_TR_L(d_L + \alpha) \biggr ] \\
& \ & + d_Td_Id_V(R_L - 1)(d_L + \alpha) \biggl [ d_T R_L (d_T R_L + 2d_I +2 d_V + d_L + \alpha) \\
& \ & + (d_L + \alpha)(d_I + 2d_V) +  (d_I + d_V)^2    \biggr ]\\
& > & d_Td_I d_V (R_L-1)(d_L + \alpha) \cdot ( d_T R_L + d_V + d_I + d_L + \alpha)^2\\
& = & A_4A_1^2 
\end{eqnarray*}
With this, all of the criteria have been satisfied and $E_I$ is stable if $R_L > 1$.  Conversely, if $R_L < 1 $, then the Jacobian possesses at least one positive eigenvalue, and the endemic state is unstable.  Finally, the local behavior of the system for $R=1$ is implied by the result of Theorem \ref{T3}.
\end{pf}

Lastly, we include a proof of the previously stated global stability theorem.
\begin{pf}[Theorem \ref{T3}]
As in \cite{Korob} for the case of (\ref{3CM}), we will prove the global stability using a Lyapunov function.
We will denote the non-infective equilibrium by $(T^0, 0,0,0)$. First, note that the quantity $T(t) - T^0 - T^0 \ln \left ( \frac{T(t)}{T^0} \right )$ vanishes when evaluated at $T(t) = T^0$ and is nonnegative as long as $T(t) > 0$ by a simple application of Taylor's Theorem.  Next, define
\begin{align*}
U(t) & = \left ( (1-p)d_L + \alpha \right ) \left [ T(t) - T^0 - T^0 \ln \left ( \frac{T(t)}{T^0} \right ) \right ] \\
 & + \left (d_L + \alpha \right ) \left [ I(t) + \frac{1}{N} V(t) \right ]  + \alpha  L(t).
\end{align*}
Notice that $U$ is nonnegative, and $U$ is identically zero if and only if it is evaluated at the non-infective equilibrium point.  We compute the derivative along trajectories and find
\vspace{-0.1in}
\begin{eqnarray*}
\frac{dU}{dt} & = & \left ( (1-p)d_L + \alpha \right ) \left ( 1 - \frac{T^0}{T}  \right ) \left [ \lambda - d_T T - kTV \right ] \\
& \ & + \left (d_L + \alpha \right )  \biggl [ (1-p) kTV + \alpha L- d_I I +  \frac{1}{N} \left ( N d_I I - d_V  V \right )  \biggr ] \\
&\ & + \alpha \left [ pkTV - (\alpha + d_L) L\right ]\\
\end{eqnarray*}
The $I, L$, and $TV$ terms all cancel and after using the definition of $T^0$, we are left with
\begin{eqnarray*}
\frac{dU}{dt}
& = & \left ( (1-p)d_L + \alpha \right ) (\lambda - d_T T ) \left ( 1 - \frac{\lambda}{d_T T} \right )\\
& \ & + \biggl [\left ( (1-p)d_L + \alpha \right )kT^* - \left (d_L + \alpha \right ) \frac{d_V}{N} \biggr ]  V  \\
& = & - \frac{ (1-p)d_L + \alpha}{d_T T} (\lambda - d_T T )^2+ \frac{(d_L  + \alpha) d_V}{N} (R_L - 1) V.
\end{eqnarray*}
Thus, under the assumption that $R_L \leq 1$, we see that $\frac{dU}{dt} \leq 0$ for all positive values of $T, I, L$, and $V$, and the global asymptotic stability follows by LaSalle's Invariance Principle.\\

Turning to the endemic equilibrium, none of the end values are zero, so we denote this steady state by $(T^*, I^*,L^*,V^*)$ and define 
\begin{align*}
U(t) & = \left ( (1-p)d_L + \alpha \right ) \left [ T(t) - T^* - T^* \ln \left ( \frac{T(t)}{T^*} \right ) \right ] \\
 & + \left (d_L + \alpha \right ) \left [ I(t) - I^* - I^* \ln \left ( \frac{I(t)}{I^*} \right )  + \frac{1}{N} \left ( V(t) - V^* - V^* \ln \left ( \frac{V(t)}{V^*} \right ) \right ) \right ]\\ 
& + \alpha  \left [ L(t) - L^* - L^* \ln \left ( \frac{L(t)}{L^*}\right ) \right ].
\end{align*}
As before, this function is nonnegative and identically zero only when evaluated at the endemic equilibrium.  Computing the derivative along trajectories yields
\begin{eqnarray*}
\frac{dU}{dt} & = & \left ( (1-p)d_L + \alpha \right ) \left ( 1 - \frac{T^*}{T}  \right ) \left [ \lambda - d_T T - kTV \right ] \\
& \ & + \left (d_L + \alpha \right ) \biggl [ \left ( 1 - \frac{I^*}{I}  \right ) \left ( (1-p) kTV + \alpha L- d_I I  \right ) \\
& \ & + \frac{1}{N} \left ( 1 - \frac{V^*}{V}  \right ) \left (N d_I I - d_V V \right )  \biggr ] + \alpha \left ( 1 - \frac{L^*}{L}  \right ) \left [ pkTV - (\alpha + d_L) L\right ]\\
& = &  \left ( (1-p)d_L + \alpha \right )  \left [ \lambda - d_T T - kTV \right ] \\
& \ & + \left (d_L + \alpha \right ) \biggl [ (1-p) kTV + \alpha L- d_I I + \left (d_I I - \frac{d_V}{N}V \right )  \biggr ]\\ 
& \ & +\alpha \left [ pkTV - (\alpha + d_L) L\right ] - \left ( (1-p)d_L + \alpha \right )  \left [ \frac{\lambda T^*}{T} - d_T T^* - kT^*V \right ] \\
& \ & - (d_L + \alpha ) \biggl [  \frac{(1-p) kTVI^*}{I} + \frac{\alpha LI^*}{I} - d_I I^*+ \frac{d_I I V^*}{V} - \frac{d_V V^*}{N} \biggr ]\\
& \ & + \alpha \left [ \frac{pkTVL^*}{L} - (\alpha + d_L) L^*\right ].
\end{eqnarray*}
Nicely, the $I, L, V$, and $TV$ terms all vanish and what remains is
\begin{align*}
\frac{dU}{dt} & =  \left ( (1-p)d_L + \alpha \right )  \left [ \lambda - d_T T + d_T T^* - \frac{\lambda T^*}{T} \right ] \\
& +  (d_L + \alpha ) \biggl [  -(1-p) k\frac{TVI^*}{I} -\alpha \frac{ LI^*}{I} + d_I I^*- d_I \frac{ I V^*}{V} + d_V\frac{ V^*}{N}\\
& + \alpha L^* - \frac{\alpha pk}{d_L + \alpha} \frac{TVL^*}{L} \biggr ]\\
& =: I + II.
\end{align*}
For $I$, we factor out a $d_T T^*$ term and use the form of $T^*$ to find
\begin{align*}
I & =  (d_L + \alpha )d_T T^* \left [R_L + 1 - \frac{T}{T^*} - R_L \frac{T^*}{T} \right ]\\
& = (d_L + \alpha )d_T T^* \left [2 - \frac{T}{T^*} - \frac{T^*}{T} + (R_L - 1) \left ( 1- \frac{T^*}{T} \right )\right ]\\
& = (d_L + \alpha )d_T T^* \left [2 - \frac{T}{T^*} - \frac{T^*}{T}\right ] + (d_L + \alpha )d_T T^*(R_L - 1) \left ( 1- \frac{T^*}{T} \right )
\end{align*}
For $II$, we factor an $L^*$ term and use the identities $$ T^*V^* = \frac{d_L + \alpha}{kp} L^* \quad \mbox{ and } \quad N d_I I^* = d_V V^*$$ to find
\begin{align*}
II & =  (d_L + \alpha )L^* \biggl [ \alpha + \frac{d_I I^*}{L^*} + \frac{d_V V^*}{N L^*}
 - (1-p)k \frac{TVI^*}{L^*I} - \frac{d_I I^*}{L^*} \frac{IV^*}{I^*V}\\ 
& \ - \frac{\alpha pk}{d_L + \alpha} \frac{TV}{L} - \alpha\frac{L I^*}{L^*I} \biggr ]\\
& = (d_L + \alpha )L^* \biggl [ \alpha + \frac{2((1-p)d_L + \alpha)}{p}
- \frac{(1-p)(d_L+\alpha)}{p} \frac{TVI^*}{T^*V^*I} \\
& - \frac{(1-p)d_L + \alpha}{p} \frac{IV^*}{I^*V}
- \alpha \frac{TVL^*}{T^*V^*L} - \alpha\frac{L I^*}{L^*I} \biggr ]\\
& = \frac{(d_L + \alpha )L^*}{p} \biggl [ ((1-p)d_L + \alpha) \left (2 - \frac{IV^*}{I^*V} \right) - (1-p)(d_L+\alpha) \frac{TVI^*}{T^*V^*I} \\
& \ + \alpha p \left ( 1 - \frac{TVL^*}{T^*V^*L} - \frac{L I^*}{L^*I} \right ) \biggr ]
\end{align*}
Thus, combining the rearrangements of $I$ and $II$, we find
\begin{align*}
\frac{dU}{dt} & = (d_L + \alpha )d_T T^* \left [2 - \frac{T}{T^*} - \frac{T^*}{T}\right ] + (d_L + \alpha )d_T T^*(R_L - 1) \left ( 1- \frac{T^*}{T} \right )\\
& \ + \frac{(d_L + \alpha )L^*}{p} \biggl [ ((1-p)d_L + \alpha) \left (2 - \frac{IV^*}{I^*V} \right) - (1-p)(d_L+\alpha) \frac{TVI^*}{T^*V^*I} \\
& \ + \alpha p \left ( 1 - \frac{TVL^*}{T^*V^*L} - \frac{L I^*}{L^*I} \right ) \biggr ]
\end{align*}
The second term simplifies to combine with those in the third term since
$$(d_L + \alpha )d_T T^*(R_L - 1)  = \frac{(d_L + \alpha )L^*((1-p)d_L + \alpha)}{p}$$
and therefore the expression becomes
\begin{align*}
\frac{dU}{dt} & = (d_L + \alpha )d_T T^* \left [2 - \frac{T}{T^*} - \frac{T^*}{T}\right ]\\
& \  + \frac{(d_L + \alpha )L^*}{p} \biggl [ ((1-p)d_L + \alpha) \left (3 - \frac{T^*}{T} - \frac{IV^*}{I^*V} \right) - (1-p)(d_L+\alpha) \frac{TVI^*}{T^*V^*I} \\
& \ + \alpha p \left ( 1 - \frac{TVL^*}{T^*V^*L} - \frac{L I^*}{L^*I} \right ) \biggr ]
\end{align*}
Since $(1-p)(d_L+\alpha) = (1-p)d_L+\alpha - \alpha p$, we add and subtract $\alpha p$ within the first term of  the second line and place the extra components in the terms on the third line to arrive at
\begin{align*}
\frac{dU}{dt}  & =  (d_L + \alpha )d_T T^* \left [2 - \frac{T}{T^*} - \frac{T^*}{T}\right ]\\
& \ + \frac{(d_L + \alpha )L^*}{p} \biggl [ (1-p)(d_L + \alpha) \left (3 - \frac{T^*}{T} - \frac{TVI^*}{T^*V^*I} - \frac{IV^*}{I^*V} \right) \\
& \ + \alpha p \left ( 4 - \frac{T^*}{T} - \frac{TVL^*}{T^*V^*L} - \frac{L I^*}{L^*I} - \frac{IV^*}{I^*V} \right ) \biggr ]
\end{align*}
Finally, each of the resulting terms above are nonpositive because the arithmetic mean is greater than the geometric mean, or more specifically, 
$$ \begin{gathered}
\frac{1}{2} \left ( \frac{T}{T^*} + \frac{T^*}{T} \right ) \geq \sqrt{\frac{T}{T^*} \cdot \frac{T^*}{T} } = 1\\
\frac{1}{3} \left ( \frac{T^*}{T} + \frac{TVI^*}{T^*V^*I} + \frac{IV^*}{I^*V} \right ) \geq \sqrt[3]{\frac{T^*}{T} \cdot \frac{TVI^*}{T^*V^*I} \cdot \frac{IV^*}{I^*V}  } = 1\\
\frac{1}{4} \left ( \frac{T^*}{T} + \frac{TVL^*}{T^*V^*L} + \frac{L I^*}{L^*I} + \frac{IV^*}{I^*V}  \right ) \geq \sqrt[4]{\frac{T^*}{T} \cdot \frac{TVL^*}{T^*V^*L} \cdot \frac{L I^*}{L^*I} \cdot \frac{IV^*}{I^*V}} = 1.
\end{gathered}$$
Hence, $\frac{dU}{dt} \leq 0$ for all positive values of $T, I, L$, and $V$. As in the non-infective case, the conclusion then follows directly from LaSalle's Invariance Principle.
\end{pf}

\section{Discussion}
In order to realistically describe and predict the effects of latent HIV infection, models of
HIV-1 dynamics and the associated mathematical tools must be capable of 
explaining the rich set of dynamics inherent within their formulation.
We have explored the steady states and asymptotic behavior of the basic three-component
model and its well-known variant which includes the effects of latent infection.
A rigorous analysis of the large time behavior of these systems displays
a reduction in the basic reproduction number due to the appearance of the latently infected T-cell population, and at 
first glance seems contradictory to the known difficulties of eradicating the latent reservoir with
antiretroviral therapy.  After undertaking a more detailed analysis here, we find that
even though the inclusion of latent T-cells allows for a wider range of parameter values to 
induce viral eradication as $t \to \infty$, the rate at which this decay occurs under ART is
retarded so significantly that, in the majority of cases, the decay could only occur outside
time periods of biological relevance.  This analysis highlights two major points.  
First and foremost, the latent cell population drastically extends the lifespan of infection.  
This can be seen from the rates of decay displayed within Section $4$ by the functions $P$ and $Q$.
However, since this property cannot be detected at the level of the basic reproduction number, a second major point becomes clear.
The standard tools of computing equilibrium states and the differing conditions under which a
system may tend to these states as $t \to \infty$ is clearly insufficient to describe, or more
importantly predict, the realistic dynamics that these equations model.  Hence, a more refined
analysis which investigates not only the end states, but the rate of propagation to a supposed 
equilibrium value within a specified time period,
is clearly needed to describe the propagation of HIV, at least when considering
the effects of latent infection.

Of course, our study is not all-inclusive.  
In attempting to address the question of latently infected cell reservoirs, we have
ignored other potential reservoirs of HIV, such as those occurring within the brain, testicles, and dendritic cells \cite{CF}.
The extent of viral replication in compartments other than resting CD4+ T-cells in patients receiving antiretroviral
therapy for extended periods of time has yet to be fully delineated.
One may also adapt the model to account for other viral reservoirs and incorporate the immune system response to a viral load. 
In addition, we assumed the use of antiretroviral therapy that included only PIs.
Certainly, the effects of RTIs could also be included, though the picture becomes slightly more complex, 
and the results are similar.
One can also study effects arising from a number of additional aspects including
\begin{enumerate}
\item A secondary infective population, such as macrophages \cite{Elaiw}
\item Pharmacological delays due to drug activation
\item The residual effects of decaying drug efficacy or periodic ART schedules
\item Spatial effects, such as those characterized by diffusion models and multiple compartment models 
\item Uncertainty arising from the measurement of parameter values or fluctuations across populations of individuals
in the form of random coefficients or stochastic differential equations 
\item Successive mutation of HIV virions
\end{enumerate}
That being said, the effects of the latent cell population on viral behavior have been clearly documented within
the current study, and it is greatly expected that even when additional mechanisms are incorporated within the
model, the basic reproduction number will not serve as a descriptive parameter alone since it only describes the
global asymptotic behavior of populations.
Future work must examine the aforementioned issues within the context of latent infection using the exponential decay functions, $P$ and $Q$,
as a refinement of the mathematical analysis detailing the long time behavior of the model.

In conclusion, the dynamics of models that consider latent infection are so complex, even when spatial fluctuations are ignored, that a single parameter, in this case $R_L$ or $R^\eps_L$, cannot possibly dictate the realistic behavior of the corresponding populations.  Instead, one must consider a number of factors including the time of validity inherent within the model, the average time periods underlying treatment, and the rates of decay associated with the trend to equilbrium.

\section{Acknowledgements}
This work is supported by the National Science Foundation under awards DMS-0908413 and DMS-1211667.  We also thank Prof. Mrinal Raghupathi (USNA) and ENS Peter Roemer (USN) for helpful comments and enthusiasm.




\end{document}